\documentclass{article}

\usepackage[utf8]{inputenc}
\usepackage[cmintegrals]{newtxmath}
\usepackage[T1]{fontenc}
\usepackage[a4paper, total={6in, 8in}]{geometry}
\usepackage{graphicx}
\usepackage{mathrsfs}
\usepackage{amsmath,stackrel}
\usepackage{graphicx}
\usepackage{multirow}
\usepackage{mathtools}
\usepackage{ marvosym }
\usepackage{tikz-cd}
\usepackage{pst-node}

\newcommand{\N}{\mathbb{N}}

\newcommand{\R}{\mathbb{R}}
\newcommand{\grad}{\text{grad}}

\newcommand{\biexp}{\text{biexp}}
\newcommand{\ind}{\text{Ind}}

\usepackage{bm}
\usepackage[colorlinks=true]{hyperref}
\hypersetup{urlcolor=blue, citecolor=red}
\newtheorem{theorem}{Theorem}[section]
\newtheorem{corollary}{Corollary}
\newtheorem{lemma}[theorem]{Lemma}
\newtheorem{proposition}{Proposition}

\newtheorem{definition}[theorem]{Definition}

\title{Local Minimizers for Variational Obstacle Avoidance on Riemannian manifolds}
\author{Jacob R. Goodman}
\date{}

\begin{document}

\maketitle

\begin{abstract}
This paper studies a variational obstacle avoidance problem on complete Riemannian manifolds. That is, we minimize an action functional, among a set of admissible curves, which depends on an artificial potential function used to avoid obstacles. In particular, we generalize the theory of bi-Jacobi fields and biconjugate points and present necessary and sufficient conditions for optimality. Local minimizers of the action functional are divided into two categories—called $Q$-local minimizers and $\Omega$-local minimizers—and subsequently classified, with local uniqueness results obtained in both cases.
\end{abstract}

\section{Introduction}
Energy-optimal path planning on nonlinear spaces such as Riemannian manifolds has been an active field of interest in the last decades due to its numerous applications in manufacturing, aerospace technologies, and robotics \cite{blochgupta, hussein, ring}. It is often the case that the desired paths must connect some set of knot points—interpolating positions with given velocities and potentially higher order derivatives \cite{CLACC, machado}. For such problems, the use of variationally defined curves has a rich history due to the regularity and optimal nature of the solutions. In particular, the so-called Riemannian {cubic} splines \cite{noakes} are a particularly ubiquitous choice in interpolant, which themselves are composed of Riemannian {cubic} polynomials—{the curves which minimize the total squared (covariant) acceleration among all sufficiently regular curves satisfying some boundary conditions in positions and velocities}—that are glued together \cite{boor}. Riemannian cubic polynomials carry a rich geometry with them which often parallels the theory of geodesics. This has been studied extensively in the literature (see \cite{Giambo, Sufficient2001} for a detailed account of Riemannian cubics and \cite{RiemannianPoly} for some results with higher-order Riemannian polynomials).  

It is often the case that—in addition to interpolating points—there are obstacles or regions in space which need to be avoided. In this case, a typical strategy is to augment the action functional with an artificial potential term that grows large near the obstacles and small away from them (in that sense, the minimizers are expected to avoid the obstacles) \cite{kod, chang, colombononholonomic}. This {strategy was used} for instance in \cite{BlCaCoCDC, BlCaCoIJC, point}, where necessary conditions for extrema in obstacle avoidance problems on Riemannian manifolds were derived, and applications to interpolation problems on manifolds and to energy-minimum problems on Lie groups and symmetric spaces endowed with a bi-invariant metric were studied. Similar strategies have been implemented for collision avoidance problems for multiagent systems evolving on Riemannian problems, as in \cite{CollAvoid, sh}. Existence of global minimizers and safety guarantees for the obstacle avoidance problem were studied in \cite{existence}. What is currently lacking in the literature regarding energy-optimal obstacle avoidance problems—and what this paper aims to address—is the derivation of sufficient conditions for optimality and a comprehensive study of the local minimization properties of the critical points.

In particular, the main contributions of this paper are as follows: (i) we investigate sufficient conditions for a curve to be a minimizing trajectory in the variational obstacle avoidance problem. This naturally leads to an extension of bi-Jacobi fields and biconjugate points as defined in \cite{Sufficient1995} for Riemannian cubics. Lemma \ref{lemma: suff_conds} provides sufficient conditions for a critical point of the action functional to be a local minimizer among the admissible set of curves (denote along the paper as an $\Omega$-local minimizer). Proposition \ref{prop: not_min} provides additional necessary conditions for the sufficient conditions to hold in terms of the so-called modified bi-Jacobi fields and biconjugate points. (ii) Robustness of the $\Omega$-local minimizers is shown in Proposition \ref{Prop: robust}, and a local uniqueness condition is provided in Proposition \ref{lemma: unique_omega_min}. This uniqueness condition is extended to the case of global minimizers in Corollary \ref{cor: unique_global}. (iii) A comprehensive study of a different type of local minimizers (called $Q$-local minimizers) is provided along section 3. In particular, a Morse index theorem is proven in Lemma \ref{lemma: index}, which is used in in Corollary \ref{cor: finite_conj} to show that there are a finite number of points biconjugate to any fixed point along a critical point of the action. Utilizing a local uniqueness condition obtained in \cite{Margarida} for the critical points of the action, we then obtain a different local uniqueness result on the restrictions of the critical points to sufficiently small intervals in Corollary \ref{cor: loc_uniq}. Finally, (iv) pairing this uniqueness with a global existence result obtained in \cite{existence}, we show in Proposition \ref{prop: loc_min} that the critical points of the action are exactly the $Q$-local minimizers.

\section{Riemannian Manifolds}\label{Sec: background}
Let $(Q, \left< \cdot, \cdot\right>)$ be an $n$-dimensional \textit{Riemannian
manifold}, where $Q$ is an n-dimensional smooth manifold and $\left< \cdot, \cdot \right>$ is a positive-definite symmetric covariant 2-tensor field called the \textit{Riemannian metric}. That is, to each point $q\in Q$ we assign a positive-definite inner product $\left<\cdot, \cdot\right>_q:T_qQ\times T_qQ\to\mathbb{R}$, where $T_qQ$ is the \textit{tangent space} of $Q$ at $q$ and $\left<\cdot, \cdot\right>_q$ varies smoothly with respect to $q$. The length of a tangent vector is determined by its norm, defined by
$\|v_q\|=\left<v_q,v_q\right>^{1/2}$ with $v_q\in T_qQ$. Given a scalar field $f:Q\to\mathbb{R}$, the metric allows one to define the \textit{gradient vector field} of $f$—denoted by $\grad f$—implicitly defined by the relation $df(q)X_q = \left<\grad f(q), X_q\right>_q$ for all $q \in Q, \ X_q \in T_q Q$.

A {\textit{connection}} $\nabla$ on $Q$ is a map that assigns to any two smooth vector fields $X$ and $Y$ on $Q$ a new vector field $\nabla_{X}Y$, playing a role similar to that of the directional derivative in classical real analysis. The operator
$\nabla_{X}$, which assigns to every vector field $Y$ the vector field $\nabla_{X}Y$, is called the \textit{covariant derivative (of $Y$) with respect to $X$}. A connection induces a number of important structures on $Q$, a particularly ubiquitous such structure is the \textit{curvature endomorphism}, which maps three vector fields $X$, $Y$ and $Z$ on $Q$ to the vector field $R(X,Y)Z := \nabla_{X}\nabla_{Y}Z-\nabla_{Y}\nabla_{X}Z-\nabla_{[X,Y]}Z.$
From an analytical perspective, the curvature endomorphism measures the extent to which covariant derivatives commute with one another. We further define the \textit{curvature tensor} $\text{Rm}$ on $Q$ via
$\text{Rm}(X, Y, Z, W) := \left<R(X, Y)Z, W\right>$.

Given a curve $q$ on $Q$ (parameterized by $t \in I \subset \R$ and with velocity vector field $\dot{q}$), there exists a unique operator $\frac{D}{dt}$ induced by $\nabla$ (called the covariant derivative along $q$) which assigns to every vector field $W$ along $q$ the vector field $\frac{D}{dt} W$ (also along $q$) which agrees with the covariant derivative $\nabla_{\dot{q}}\tilde{W}$ for any extension $\tilde{W}$ of $W$ to $Q$. A vector field $X$ along $q$ is said to be \textit{parallel along $q$} if $\displaystyle{\frac{DX}{dt}\equiv 0}$. For $k\in \N$, the $k$th-order covariant derivative of $W$ along $q$, denoted by $\displaystyle{\frac{D^k}{dt^k} W}$, can then be inductively defined by $\displaystyle{\frac{D^k}{dt^k} W = \frac{D}{dt} \left(\frac{D^{k-1}}{dt^{k-1}} W\right)}$.

It is well-known that the Riemannian metric induces a unique torsion-free and metric compatible connection called the \textit{Riemannian connection}, or the \textit{Levi-Civita connection}. Along the remainder of this paper, we will assume that $\nabla$ is the Riemannian connection. For additional information on connections and curvature, we refer the reader to \cite{Boothby,Milnor}.

If we assume that $Q$ is \textit{complete}, then {by the Hopf-Rinow theorem,} any two points $x$ and $y$ in $Q$ can be connected by a (not necessarily unique) minimal-length curve $\gamma_{x,y}$. In particular, $\gamma_{x,y}$ must be a \textit{geodesic}—that is, it verifies $\frac{D}{dt}\dot{\gamma}_{x, y} \equiv 0$. In this case, the Riemannian distance between $x$ and $y$ can be defined by
{$\displaystyle{d(x,y)=\int_{0}^{1}\Big{\|}\frac{d \gamma_{x,y}}{d s}(s)\Big{\|}\, ds}$}. Geodesics provide a map $\mathrm{exp}_q:T_qQ\to Q$ called the \textit{Riemannian exponential map} such that $\mathrm{exp}_q(v) = \gamma(1)$, where $\gamma$ is the unique geodesic verifying $\gamma(0) = q$ and $\dot{\gamma}(0) = v$. In particular, $\mathrm{exp}_q$ is a diffeomorphism from some star-shaped neighborhood of $0 \in T_q Q$ to a convex open neighborhood $\mathcal{B}$ of $q \in Q$. If $y \in \mathcal{B}$, we can write the Riemannian distance by means of the Riemannian exponential as $d(q,y)=\|\mbox{exp}_q^{-1}y\|.$ 

Let $Q$ be an $m$-dimensional Riemannian manifold and $a < b \in \R$. Given $\xi = (q_a, v_a), \eta = (q_b, v_b) \in TQ$, we denote the space of piece-wise smooth $C^1$ curves $\gamma: [a, b] \to Q$ satisfying $\gamma(a) = q_a, \ \gamma(b) = q_b, \ \dot{\gamma}(a) = v_a, \ \dot{\gamma}(b) = v_b$ by $\Omega_{\xi, \eta}^{a, b}$. Along the paper, we will frequently drop the subscripts and superscripts on $\Omega_{\xi, \eta}^{a,b}$ unless otherwise necessary. $\Omega$ has the structure of a smooth manifold, and the tangent space $T_x \Omega$ consists of all piece-wise smooth $C^1$ vector fields $X$ along $x$ satisfying $X(a) = X(b) = \frac{D}{dt}X(a) = \frac{D}{dt}X(b) = 0$. Such a vector field can be viewed as the variational vector field of some $\textit{variation}$ of $x$, which is a function $\alpha: (-\epsilon, \epsilon) \times [a, b] \to Q$ such that:
\begin{enumerate}
    \item $\alpha(r, \cdot) \in \Omega$ for all $r \in (-\epsilon, \epsilon)$,
    \item $\alpha(0, t) = x(t)$ for all $t \in [a, b]$,
    \item $\frac{d}{dr}\big\vert_{r=0}\alpha(r, t) = X(t)$ for all $t \in [a, b]$.
\end{enumerate}
Consider the norm $\| \cdot \|_{T_x \Omega}$ on $T_x \Omega$ given by
$$\|X\|_{T_x \Omega} = \left(\int_a^b \left[\left\|X\right\|^2 + \left\|\frac{DX}{dt}\right\|^2 + \left\|\frac{D^2 X}{dt^2}\right\|^2 \right]dt\right)^{1/2}.$$
We let $\mathring{H}^2_x$ be the completion of $T_x \Omega$ under $\| \cdot \|_{T_x \Omega}$. Considering an orthonormal basis of parallel vector fields $\left\{X_i\right\}$ along $x$ and writing $X = \xi^i X_i$ for some $\xi^i: [a, b] \to \R$, we have that
$$\|X\|_{T_x \Omega} = \left(\int_a^b \left[ \xi^i \xi^i + \dot{\xi}^i \dot{\xi}^i + \ddot{\xi}^i \ddot{\xi}^i \right]dt\right)^{1/2},$$
from which it is clear that $\mathring{H}^2_x$ can be identified with the Sobolev space $\mathring{H}^2([a,b], \R^n)$ (as discussed for instance in section 4.3 of \cite{Jost}).

\section{Necessary and Sufficient Conditions for the Variational Obstacle Avoidance Problem}\label{Sec: Necessary conditions}


Consider a complete and connected Riemannian manifold $Q$, and for some $\xi, \eta \in TQ$ and $T > 0$, let $\Omega_{\xi, \eta}^{0, T}$ be defined as in the previous section. We define the function $J: \Omega \to \R$ by:
\begin{equation}\label{J}
J(q)=\int\limits_0^T \Big{(}\frac12\Big{|}\Big{|}\frac{D }{dt}\dot{q}(t)\Big{|}\Big{|}^2 + V(q(t))\Big{)}dt.
\end{equation}

\textbf{Variational obstacle avoidance problem:} Find a curve $q\in \Omega$ minimizing the functional $J$, where $V:Q \to\mathbb{R}$ is a smooth and non-negative function called the \textit{artificial potential}.

\vspace{.2cm}

In order to minimize the functional $J$ among the set $\Omega$, we want to find curves $q\in\Omega$ such that $J(q)\leq J(\tilde{q})$ for all
admissible curves $\tilde{q}$ in a $C^1$-neighborhood of $q$. Necessary conditions can be derived by finding $q$ such that the differential of $J$ at $q$, $dJ(q)$, vanishes identically. This is clearly equivalent to $dJ(q)W = 0$ for all $W \in T_q \Omega$, which itself can be understood through variations—as discussed in the previous section. The next result from \cite{BlCaCoCDC} characterizes necessary conditions for optimality in the variational obstacle avoidance problem.


\begin{proposition}\label{th1}\cite{BlCaCoCDC} A point $q \in \Omega$ is a critical point for the functional $J$ if and only if it is a $\mathcal{C}^\infty$-curve on $[0,T]$ satisfying
\begin{equation}\label{eqq1}
    \frac{D^3\dot{q}}{dt^3}+R\Big{(}\frac{D\dot{q}}{dt},\frac{dq}{dt}\Big{)}\frac{dq}{dt} + \hbox{\grad} \, V(q(t)) = 0.
\end{equation}
\end{proposition}

We will call solutions to equation \eqref{eqq1} \textit{modified cubic polynomials}.  It is natural to ask: when are modified cubic polynomials minimizers of $J$? This can be understood both locally and globally—although there is some discrepancy in how these terms are used. For that reason, we introduce the following definition:

\begin{definition}
A curve $q \in \Omega$ is a:
\begin{enumerate}
    \item[\hbox{(i)}] \textit{Global minimizer} of $J$ on $\Omega$ iff $J(q) \le J(\tilde{q})$ for all $\tilde{q} \in \Omega$.
    \item[\hbox{(ii)}] $\Omega$-local minimizer of $J$ on $\Omega$ iff $J(q) \le J(\tilde{q})$ for all $\tilde{q}$ in some $C^1$ neighborhood of $q$ (within $\Omega$).
    \item[\hbox{(iii)}] $Q$-local minimizer of $J$ iff for any $\tau \in [0, T]$, there exists an interval $[a, b] \subset [0, T]$ containing $\tau$ such that $q\vert_{[a,b]}$ is a minimizer of $J$ on $\Omega_{\xi, \eta}^{[a, b]}$, where $\xi = (q(a), \dot{q}(a)), \ \eta = (q(b), \dot{q}(b)) \in TQ$.
\end{enumerate}
\end{definition}

It should be noted that we have slightly abused our notation in the definition of a $Q$-local minimizer. Technically, we are concerned with minimizing the integral $\int_a^b \Big{(}\frac12\Big{|}\Big{|}\frac{D }{dt}\dot{q}(t)\Big{|}\Big{|}^2 + V(q(t))\Big{)}dt$, which has different limits of integration than $J$ as defined in equation \eqref{J}. We will continue to refer to integrals of this form by $J$ throughout the paper—and in every case, the limits of integration will match that of the boundary conditions defined by the admissible set $\Omega_{\xi, \eta}^{[a, b]}$ on which $J$ is being discussed.

\subsection{$\Omega$-local minimizers of $J$}

In order to understand when a critical point of $J$ (that is, a modified cubic polynomial) is an $\Omega$-local minimizer of $J$, we utilize \textit{second variations}, akin to the classical second derivative test from calculus.
In particular, for some modified cubic polynomial $q$, and some real numbers $\epsilon_r, \epsilon_s > 0$, we consider a \textit{two-parameter variation} $\alpha: (-\epsilon_r, \epsilon_r) \times (-\epsilon_s, \epsilon_s) \times [0, T] \to Q$ of $q$ such that:
\begin{enumerate}
    \item $\alpha(r, s, \cdot) \in \Omega$ for all $(r, s) \in (-\epsilon_r, \epsilon_r) \times (-\epsilon_s, \epsilon_s)$,
    \item $\alpha(0, 0, t) = q(t)$ for all $t \in [0, T]$,
    \item $\frac{d}{dr}\big\vert_{r=0}\alpha(r, 0, \cdot) = X \in T_{q} \Omega$,
    \item $\frac{d}{ds}\big\vert_{s=0}\alpha(0, s, \cdot) = Y \in T_{q} \Omega$.
\end{enumerate}
A particularly ubiquitous choice is given by $\alpha(r, s, t) = \exp_{q(t)}(rX(t) + sY(t))$. The second variation of $J$ along a modified cubic polynomial $q$ is then calculated as $\frac{\partial^2}{\partial r\partial s}\Big\vert_{(r,s)=(0,0)} J(\alpha(r, s, \cdot)) = \frac{\partial}{\partial r}\Big\vert_{r=0}\frac{\partial}{\partial s}\Big\vert_{s=0} J(\alpha(r, s, \cdot))$. From \cite{BlCaCoCDC} Theorem 3.1, we see that

\begin{align*}
    \frac{\partial}{\partial s}J(\alpha) &= \int_0^T \left<\frac{\partial \alpha}{\partial s}, \frac{D^4\alpha}{dt^4} + R\left(\frac{D^2}{dt^2}\alpha, \frac{\partial\alpha}{dt} \right)\frac{\partial\alpha}{dt} + \grad V(\alpha) \right>dt \\
    &+ \sum_{i=1}^l \left[\left<\frac{D}{dt}\frac{\partial \alpha}{\partial s}, \frac{D^2 \alpha}{dt^2}\right> - \left<\frac{\partial \alpha}{\partial s}, \frac{D^3 \alpha}{dt^3} \right>\right]_{t_{i-1}^+}^{t_i^-},
\end{align*}
and therefore,
\begin{align}
    \frac{\partial^2}{\partial r\partial s}\Big\vert_{(r,s)=(0,0)} J(\alpha) &= \int_0^T \left<\frac{D}{\partial r}\frac{\partial \alpha}{\partial s}\Big\vert_{(r,s)=(0,0)}, \ \frac{D^4 q}{dt^4} + R\left(\frac{D^2 q}{dt^2}, \dot{q} \right)\dot{q} + \grad V(q) \right>dt\\
    &+ \int_0^T \left<Y, \frac{D}{\partial r}\left[ \frac{D^4\alpha}{dt^4} + R\left(\frac{D^2}{dt^2}\alpha, \frac{\partial\alpha}{dt} \right)\frac{\partial\alpha}{dt}\right] \right>dt + \frac{D}{\partial r}\sum_{i=1}^l \left[\left<\frac{D}{dt}\frac{\partial \alpha}{\partial s}, \frac{D^2 \alpha}{dt^2}\right> - \left<\frac{\partial \alpha}{\partial s}, \frac{D^3 \alpha}{dt^3} \right>\right]_{t_{i-1}^+}^{t_i^-} \\
    &+ \int_0^T \left<Y(t), \frac{D}{\partial r}\Big\vert_{r=0}\grad V(\alpha(r, 0, t))\right>dt.
\end{align}

The first integral vanishes identically since $q$ is a modified cubic polynomial. The second integral was calculated in \cite{Sufficient1995}, which studies sufficient conditions for optimality for Riemannian cubic polynomials (that is, with $V \equiv 0$). Hence, we need only calculate the final integral. If we define $\gamma_t(r) := \alpha(r, 0, t)$ and denote $\dot{\gamma}_t(r) = \frac{d}{dr}\gamma_t(r)$, then $\frac{D}{\partial r}\Big\vert_{r=0} \grad V(\alpha(r, 0, t)$ can be expressed as:
$$\frac{D}{\partial r}\Big\vert_{r=0} \grad V(\gamma_t(r)) = \nabla_{\dot{\gamma_t}(r)} \grad V \Big\vert_{r=0} = \nabla_{X(t)} \grad V.$$
Adding this term along with those found in \cite{Sufficient1995} Theorem 2.4, we obtain the following lemma:
\begin{lemma}\label{lemma: suff_conds}
\begin{align*}
    \frac{\partial^2}{\partial r\partial s}\Big\vert_{(r,s)=(0,0)} J(\alpha) &= \int_0^T \left< Y(t), \frac{D^4}{dt^4} X  + F(X(t), \dot{q}(t)) + \nabla_{X(t)}\grad V(q(t))\right>dt \\ 
    &- \sum_{i=1}^l \left<\frac{D}{dt}Y(t_i), \frac{D^2}{dt^2} X(t_i^+) - \frac{D^2}{dt^2}X(t_{i}^-) \right> + \sum_{i=1}^l \left<Y(t_i), \frac{D^3}{dt^3} X(t_i^+) - \frac{D^3}{dt^3}X(t_{i}^-) \right> 
\end{align*}
Where 
\begin{align*}
    F(X, Y) &= (\nabla^2_Y R)(X, Y)Y + (\nabla_X R)(\nabla_Y Y, Y)Y + R(R(X, Y)Y, Y)Y + R(X, \nabla^2_Y Y)Y \\
    &+ 2\left[(\nabla_Y R)(\nabla_Y X, Y)Y + (\nabla_Y R)(X, \nabla_Y Y)Y + R(\nabla^2_Y X, Y)Y \right] \\
    &+ 3\left[ (\nabla_Y R)(X, Y)\nabla_Y Y + R(X, Y) \nabla_Y^2 Y + R(X, \nabla_Y Y)\nabla_Y Y \right] + 4R(\nabla_Y X, Y)\nabla_Y Y
\end{align*}
\end{lemma}

We further define the bilinear form $I: T_q \Omega \times T_q \Omega \to \R$ called the \textit{index form} by
\begin{align*}
    I(X, Y) = \int_0^T \left[\left<\frac{D^2 X}{dt^2}, \frac{D^2 Y}{dt^2} \right> + \left<Y, F(X, \dot{q}) + \nabla_X \grad V \right>\right]dt,
\end{align*}
from which it can be easily observed by integrating by parts twice that $\frac{\partial^2}{\partial r\partial s}\Big\vert_{(r,s)=(0,0)} J(\alpha) = I(X, Y)$. However, unlike the case of Riemannian cubics and geodesics where the corresponding index form is symmetric, the artificial potential introduces an asymmetry. Namely, we can decompose $I$ as $I(X, Y) = I_c(X, Y) + P_+(X, Y) + P_-(X,Y)$, where 
\begin{itemize}
    \item $I_c(X, Y) = \int_0^T \left[\left<\frac{D^2 X}{dt^2}, \frac{D^2 Y}{dt^2} \right> + \left<Y, F(X, \dot{q}(t)) \right>\right]dt$ is the symmetric bilinear form obtained for Riemannian cubics in \cite{Sufficient1995, Sufficient2001}.
    \item $P_+(X, Y) = \frac12\int_0^T \left[ \left<Y, \nabla_X \grad V\right> + \left<X, \nabla_Y \grad V\right>\right]dt$ is a symmetric bilinear form describing the symmetric contribution of the artificial potential.
    \item $P_-(X, Y) = \frac12\int_0^T \left[\left<Y, \nabla_X \grad V\right> - \left<X, \nabla_Y \grad V\right>\right]dt$ is an anti-symmetric bilinear form describing the anti-symmetric contribution of the artificial potential.
\end{itemize}  

A sufficient condition for a modified cubic polynomial $q$ to be an $\Omega$-local minimizer of $J$ is that $\frac{\partial^2}{\partial r\partial s}\Big\vert_{(r,s)=(0,0)} J(\alpha) \ge 0$ for all $2$-parameter variations $\alpha$ along $q$. Or, equivalently, that the index form is positive semi-definite for all variational vector fields along $q$. Of course, such a calculation is daunting in practice, so we now turn our attention to understanding when $q$ is \textit{not} an $\Omega$-local minimizer by studying the vector fields belonging to the kernel of $I$. In particular, this yields the analogue of Jacobi fields for modified cubic polynomials.

\begin{lemma}\label{lemma: Jacobi}
A vector field $X \in T_q \Omega$ belongs to the kernel of $I$ if and only if $X$ is smooth and $\frac{D^4}{dt^4} X + F(X, \dot{q}) + \nabla_X \grad V(q) \equiv 0$.
\end{lemma}

\textit{Proof:} Clearly, if $X$ is smooth satisfies $\frac{D^4}{dt^4} X + F(X, \dot{q}) + \nabla_X \grad V(q) \equiv 0$, then $\frac{\partial^2}{\partial r\partial s}\Big\vert_{(r,s)=(0,0)} J(\alpha) = I(X, Y) = 0$ for any $Y \in T_q \Omega$. On the other hand, if $X$ belongs to the kernel of $I$, then choose $Y(t) = f(t)(\frac{D^4}{dt^4} X + F(X, \dot{q}) + \nabla_X \grad V(q))$, where $f \in C^{\infty}([0, T])$ is non-negative and satisfies $f(t_i) = \frac{d}{dt}f(t_i) = 0$. Then we have
$\int_0^T f(t) \left\|Y(t) \right\|^2 = 0$, which implies $Y = 0$ for all $t \in (t_i, t_{i+1})$ for $1 \le i \le l-1$. To see that $X$ is smooth, we first set $Y(t) = f(t)Z(t)$, where $f(t_i) = 0, \ \frac{d}{dt} f(t_i) > 0$ and $Z(t_i) = \frac{D^2}{dt^2}X(t_i^+) - \frac{D^2}{dt^2}X(t_i^-)$. Then, we obtain $\frac{d}{dt}f(t_i)\|Z(t_i)\|^2 = 0 \implies \frac{D^2}{dt^2}X(t_i^+) = \frac{D^2}{dt^2}X(t_i^-)$. Next, we set $Y(t) = f(t)Z(t)$ with $f(t_i) > 0$ and $Z(t_i) = \frac{D^3}{dt^3}X(t_i^+) - \frac{D^3}{dt^3}X(t_i^-)$, from which we find that $f(t_i)\|Z(t_i)\|^2 = 0 \implies \frac{D^3}{dt^3}X(t_i^+) = \frac{D^3}{dt^3}X(t_i^-)$. Hence $X$ is smooth. $\hfill\square$

This motivates the following definition:
\begin{definition}
A vector field $X$ along a modified cubic $q$ satisfying $\frac{D^4}{dt^4} X + F(X, \dot{q}) + \nabla_X \grad V(q) \equiv 0$ on $[0, T]$ is called a modified bi-Jacobi field.
\end{definition}
Observe that in the case where $V \equiv 0$, the definition of a modified bi-Jacobi Field coincides with that of a bi-Jacobi field, as defined in \cite{Sufficient2001}. Moreover, note that the equation describing the modified bi-Jacobi fields is linear in $X$, so that (since $V$ is smooth) the modified bi-Jacobi fields are smooth and the existence and uniqueness of solutions on $[0, T]$ given initial values $X(0), \ \frac{D}{dt}X(0), \ \frac{D^2}{dt^2}X(0), \ \frac{D^3}{dt^3}X(0)$ follows immediately (say, by moving to coordinate charts). In particular, the set of modified bi-Jacobi fields along a modified cubic polynomial $q$ forms a $4n$-dimensional vector space. 

\begin{definition}
Two points $t = t_1, t_2 \in [0, T]$ are said to be biconjugate along a modified cubic $q$ if there exists a non-zero modified bi-Jacobi field $X$ such that
\begin{align*}
    X(t_1) = X(t_2) = 0, \quad \text{ and }\qquad \frac{D}{dt}X(t_1) = \frac{D}{dt}X(t_2) = 0.
\end{align*}
\end{definition}

\noindent Analogous to the case of geodesics and conjugate points (\cite{Jost}, Theorem 4.3.1), or Riemannian cubic polynomaials and biconjugate points (\cite{Sufficient2001}, Theorem 7.2), we now show that modified cubic polynomials do not minimize past their biconjugate points.

\begin{proposition}\label{prop: not_min}
Suppose that $q \in \Omega$ is a modified cubic polynomial and $0 \le t_1 < t_2 < T$ are biconjugate. Then there exists a vector field $X \in T_q \Omega$ such that $I(X, X) < 0$. In particular, $q$ is not a minimizer of $J$ on $\Omega$.
\end{proposition}

\textit{Proof:}
Since $t_1, t_2$ are biconjugate, there exists a modified bi-Jacobi field $U$ such that $U(t_1) = U(t_2) = \frac{D}{dt}U(t_1) =\frac{D}{dt}U(t_2) = 0$. We then consider the vector field 
$X(t) = \begin{cases} U(t) &t\in[t_1, t_2] \\
0, &\text{otherwise}
\end{cases}$. It is clear that $X \in T_q \Omega$ is smooth except possibly at the points $t = t_1, t_2$, and that at least one of $\frac{D^2}{dt^2} U(t_1), \frac{D^3}{dt^3} U(t_1)$ are non-zero and at least one of $\frac{D^2}{dt^2} U(t_2), \frac{D^3}{dt^3} U(t_2)$ are non-zero (otherwise, $U \equiv 0$). Now, we define two smooth vector fields $Z, W$ along $\dot{q}$ such that
\begin{align*}
    Z(t_1) &= -\frac{D^3}{dt^3} U(t_1), \qquad Z(t_2) = -\frac{D^3}{dt^3} U(t_2), \qquad \frac{D}{dt}Z(t_1) = \frac{D}{dt}Z(t_2) = 0, \\
    W(t_1) &= W(t_2) = 0,  \qquad \frac{D}{dt}W(t_1) = \frac{D^2}{dt^2} U(t_1), \qquad  \frac{D}{dt}W(t_2) = \frac{D^2}{dt^2} U(t_2),
\end{align*}
and define $Y(t) = \phi(t)Z(t) + \psi(t)W(t)$, where $\phi, \psi \in C^{\infty}([0, T])$ such that
\begin{align*}
    \phi(t_1) &= \phi(t_2) = 1, \qquad \dot{\phi}(t_1) = \dot{\phi(t_2)} = 0, \qquad 0 < \phi \le 1 \quad \forall t \in \text{supp}(\phi) \subseteq (t_1 - \delta, t_1 + \delta) \cup (t_2 - \delta, t_2 + \delta),\\
    \psi(t_1) &= \psi(t_2) = 0, \qquad \dot{\psi}(t_1) = \dot{\psi(t_2)} = 1, \qquad 0 < \psi \le 1 \quad \forall t \in \text{supp}(\psi) \subseteq (t_1 - \delta, t_1 + \delta) \cup (t_2 - \delta, t_2 + \delta),
\end{align*}
Then, in particular $Y$ satisfies:
\begin{align*}
    Y(t_1) = -\frac{D^3}{dt^3} U(t_1), \qquad \frac{D}{dt}Y(t_1) = \frac{D^2}{dt^2} U(t_1), \qquad Y(t_2) = -\frac{D^3}{dt^3} U(t_2), \qquad \frac{D}{dt}Y(t_2) = \frac{D^2}{dt^2} U(t_2),
\end{align*}
with $\text{supp}(Y) \subset (t_1 - \delta, t_1 + \delta) \cup (t_2 - \delta, t_2 + \delta)$. For $\epsilon > 0$, we consider the vector field $U_{\epsilon} = X + \epsilon Y$. It is clear that $U_{\epsilon} \in T_q \Omega$, and by the bilinearity of $I$, we have $I(U_{\epsilon},U_{\epsilon}) = I(X, X) + \epsilon( I(X, Y) + I(Y, X)) + \epsilon^2 I(Y, Y)$. Since $X$ is smooth except at $t = t_1, t_2$, where it vanishes (along with its covariant derivative) and it satisfies the bi-Jacobi equation on $[0, t_1], [t_1, t_2], [t_2, T]$, it is clear that $I(X, X) = 0.$ Moreover, we have that $I(X, Y) = -\left\|Y(t_1)\right\|^2 - \left\|\frac{D}{dt}Y(t_1)\right\|^2 -\left\|Y(t_2)\right\|^2 - \left\|\frac{D}{dt}Y(t_2)\right\|^2 < 0$. Now observe that 
\begin{align*}
    I(Y, X) &= I_c(Y, X) + P_+(Y, X) + P_-(Y, X) \\
    &= I_c(X, Y) + P_+(X, Y) - P_-(X, Y) \\
    &= I(X, Y) - 2P_-(X, Y)
\end{align*}

Consider the 2-form $G(X, Y) = \left<Y, \nabla_X \grad V\right> - \left<X, \nabla_Y \grad V\right>$. Then we have:
$$|G(X, Y)| \le \phi |G(X, Z)| + \psi |G(X, W)| \le (\phi + \psi)(\|Z\| \|\nabla_W \grad V\| + \|W\| \|\nabla_Z \grad V\|),$$
hence \begin{align*}|P_-(X, Y)| \le \int_0^T |G(X, Y)|dt &\le 2\delta \max_{t \in [0, T]}\left\{\|Z(t)\| \|\nabla_W(t) \grad V(q(t))\| + \|W(t)\| \|\nabla_Z(t) \grad V(q(t))\| \right\} \\
&:= \delta C.
\end{align*}
It follows that $I(Y, X) \le I(X, Y) + 2\delta C$. Therefore, $I(U_{\epsilon}, U_{\epsilon}) \le 2\epsilon(I(X, Y) + \delta C) + \epsilon^2 I(Y, Y)$. If we choose $\delta$ so that $I(X, Y) + \delta C < 0$, then it is clear that this quantity is negative for sufficiently small $\epsilon$. $\hfill\square$

One might wish to understand the robustness of being an $\Omega$-local minimizer. That is, if $q$ is an $\Omega$-local minimizer, is it also true that the restriction of $q$ to subset of $[0, T]$ is an $\Omega$-local minimizer on its corresponding admissible set? It turns out that the answer is yes, which is summarized in the following proposition.

\begin{proposition}\label{Prop: robust}
Suppose that $q$ is an $\Omega$-local minimizer of $J$ and let $[a, b] \subset [0, T]$. Then the curve $q \vert_{[a,b]}$ is an $\Omega_{\xi, \eta}^{a, b}$-local minimizer, where $\xi = (q(a), \dot{q}(a)), \ \eta = (q(b), \dot{q}(b)) \in TQ$. 
\end{proposition}

\noindent\textit{Proof:} Since $q$ is an $\Omega$-local minimizer, there exists a $C^1$ neighborhood $B$ of $q$ contained in $\Omega$ such that $J(q) \le J(\tilde{q})$ for any $\tilde{q} \in B$. Let $B^\ast$ be the set of curves $q^\ast \in \Omega_{\xi, \eta}^{a, b}$ such that the curve $\tilde{q}(t) = \begin{cases} q(t), & t\in[0, a) \cup [b, T] \\ 
q^\ast(t), & t \in [a, b]\end{cases}$ is contained in $B$. Clearly, $B^\ast$ is a non-empty $C^1$ neighborhood of $q^\ast$ contained in $\Omega_{\xi, \eta}^{a, b}$. Fix a $q^\ast \in B^\ast$, and the corresponding $\tilde{q} \in B$. Then by the additivity of integration and the fact that $J(q) \le J(\tilde{q})$, it follows that
\begin{align*}
    J(q \vert_{[0, a]}) + J(q \vert_{[a, b]}) + J(q \vert_{[b, T]}) &\le J(\tilde{q}\vert_{[0, a]}) + J(\tilde{q}\vert_{[a, b]}) + J(\tilde{q}\vert_{[b, T}), \\
    J(q \vert_{[0, a]}) + J(q \vert_{[a, b]}) + J(q \vert_{[b, T]}) &\le J(q\vert_{[0, a]}) + J(q^\ast\vert_{[a, b]}) + J(q\vert_{[b, T}), \\
     J(q\vert_{[a, b]}) &\le  J(q^\ast\vert_{[a,b]}).
\end{align*} 
Hence, $q\vert_{[a, b]}$ is an $\Omega_{\xi, \eta}^{a, b}$-local minimizer of $J$. $\hfill \square$

Denote by $Tq: [0, T] \to TQ$ the \textit{tangent lift} of the curve $q \in \Omega$. That is, the curve lying in the tangent bundle $TQ$ whose local coordinate expression is given by $(q^i, \dot{q}^i)$. 

\begin{proposition}
\label{lemma: unique_omega_min}
Suppose that $q$ is an $\Omega$-local minimizer of $J$, and let $B$ be a $C^1$ neighborhood of $q$ contained in $\Omega$ such that $J(q) \le J(\tilde{q})$ for any $\tilde{q} \in B$. If $q_0 \in B$ satisfies $J(q_0) = J(q)$ and $Tq_0(\tau) = Tq(\tau)$ for some $\tau \in (0, T)$, then $q_0 \equiv q$ on $[0, T]$. 
\end{proposition}

\noindent\textit{Proof:} Suppose that $q_0 \in B$ is such that $Tq(\tau) = Tq_0(\tau)$ and $J(q) = J(q_0)$, and consider the curve $\tilde{q}(t) = \begin{cases} q(t), & t\in[0, \tau) \\ 
q_0(t), & t \in [\tau, T]\end{cases}$. Clearly, $\tilde{q}$ is contained in $B$ and $J(q) = J(\tilde{q})$, so that $\tilde{q}$ is an $\Omega$-local minimizer (on $B$, for instance) and hence is a critical point of $J$ on $\Omega$. But by Proposition \ref{th1}, this implies that $\tilde{q}$ is smooth, and therefore
$$\begin{cases}q(\tau) = q_0(\tau) \\
\dot{q}(\tau) = \dot{q}_0(\tau) \\
\frac{D}{dt}\dot{q}(\tau) = \frac{D}{dt}\dot{q}_0(\tau) \\
\frac{D^2}{dt^2}\dot{q}(\tau) = \frac{D^2}{dt^2}\dot{q}_0(\tau).\end{cases}$$
By the uniqueness of solutions to equation \eqref{eqq1} with the above initial conditions, we have that $q \vert_{[\tau, T]} \equiv q_0 \vert_{[\tau, T}]$. Repeating this argument backwards in time, we reach the desired conclusion. $\hfill \square$

Lemma \ref{lemma: unique_omega_min} can be interpreted as saying that the tangent lifts of sufficiently nearby $\Omega$-local minimizers cannot intersect away from the boundary points. Alternatively, this means that $\Omega$-local minimizers are (locally) uniquely determined by a single internal position and velocity. If we were to consider this Proposition \ref{lemma: unique_omega_min} in the case that $q$ is a \textit{global} minimizer, we could of course extend $B$ to be the entire space $\Omega$. Together with Proposition \ref{Prop: robust}, this leads to the following corollary.

\begin{corollary}\label{cor: unique_global}
Suppose that $q \in \Omega_{\xi, \eta}^{a, b}$ is a global minimizer of $J$, and let $[a^\ast, b^\ast] \subset [a, b]$ be a proper subset. Then, for $\mu = (q(a^\ast), \dot{q}(a^\ast)), \ \nu = (q(b^\ast), \dot{q}(b^\ast)) \in T Q$, the curve $q \vert_{[a^\ast, b^\ast]}$ is the unique minimizer of $J$ on $\Omega_{\mu, \nu}^{a^\ast, b^\ast}$
\end{corollary}

\subsection{$Q$-local minimizers of $J$}

We now wish to understand when the critical points are $Q$-local minimizers of $J$. In Proposition \ref{prop: loc_min}, we will show that every modified cubic polynomial is a $Q$-local minimizer of $J$—just as every geodesic is a $Q$-local minimizer of the length functional. In particular, we will do this by first establishing a local uniqueness result for modified cubic polynomials, based off of a similar result for Riemannian cubic polynomials provided in \cite{Margarida}. 

Consider the $3$-tangent bundle $(TQ)^3 := \cup_{p \in Q} (T_p Q)^3$. It is easily seen by (for example by considering local coordinate charts and applying standard ODE theory) that a modified cubic polynomial is uniquely determined by initial conditions
$$q(0) = p, \quad \dot{q}(0) = v, \quad \frac{D}{dt} \dot{q}(0) = y, \quad \frac{D^2}{dt^2}\dot{q}(0) = z,$$
for $(p, v, y, z) \in (TQ)^3$. That is, there exists a neighborhood $B \subset (TQ)^3$ containing $(p, v, y, z)$ and a $\delta > 0$ such that for each $(p_0, v_0, y_0, z_0) \in B$, there exists a unique modified cubic polynomial $q: (-\delta, \delta) \to Q$ with the initial conditions prescribed by $(p_0, v_0, y_0, z_0)$. If we consider a set $C \subset (T_p Q)^2$ defined such that $\{(p, v)\} \times C \subset B$, we may then consider the \textit{bi-exponential map} $\biexp_{(p,v)}^t: C \to TQ$ defined by $\biexp_{(p,v)}^t(y, z) = (q(t), \dot{q}(t))$, where $q$ is the unique modified cubic polynomial satisfying the initial conditions prescribed by $(p, v, y, z) \in (TQ)^3$, and $t \in (0, \delta)$.

The bi-exponential map serves as a connection between the initial value problem and the boundary value problem. As we will see below, the differential of the bi-exponential map determines when we have local uniqueness of solutions to equation \eqref{eqq1} (that is, when modified cubic polynomials are locally unique). Furthermore, there is an intimate connection between biconjugate points and the differential of the bi-exponential map, described by Lemma \ref{lemma: conjugate}. Lemma \ref{lemma: unique_exp} and Lemma \ref{lemma: conjugate} appear in \cite{Margarida} in the case that $V \equiv 0$, however the proofs follow almost identically in both cases.\\

\begin{lemma}\label{lemma: unique_exp}
Let $(p, v, y, z) \in (TQ)^3$ and $\biexp_{(p,v)}^t$ be defined in a neighborhood $C$ of $(y, z) \in (T_p Q)^2$ for $t \in (0, \delta)$, with $\delta > 0$. If $\biexp_{(p,v)}^{\tau}$ is not critical at $(y, z)$, then there exists a neighborhood $W_1$ of $(y, z)$, a neighborhood $W_2$ of $\biexp_{(p,v)}^{\tau}(y, z)$, and a neighborhood $V$ of $(p, v, y, z)$ such that for each $(p_1, v_1) \in W_1$ and $(p_2, v_2) \in W_2$, there exists a unique modified cubic polynomial $q$ satisfying
$$q(0) = p_1, \quad \dot{q}(0) = v_1, \quad q(\tau) = p_2, \quad \dot{q}(\tau) = v_2,$$
and $(q(0), \ \dot{q}(0), \ \frac{D}{dt} \dot{q}(0), \ \frac{D^2}{dt^2} \dot{q}(0)) \in V$.
\end{lemma}

\begin{lemma}\label{lemma: conjugate}
$biexp_{(p,v)}^{\tau}$ is not critical at $(y, z) \in (T_p Q)^2$ if and only if the points $t=0, \tau$ are not biconjugate along the modified cubic polynomial $t \mapsto \pi \circ \biexp_{(p,v)}^{t}(y, z)$.
\end{lemma}

We now wish to establish that there are a finite number of points biconjugate to any $\tau \in [0, T].$ For this, we turn to some techniques of global analysis, as discussed for instance in \cite{Jost}. For that, we consider the index form as a quadratic form on the Hilbert space $\mathring{H}^2_q$, as discussed in Section \ref{Sec: background}. That is, we consider $I: \mathring{H}^2_q \times \mathring{H}^2_q \to \R$ defined by:
\begin{align*}
    I(X, Y) = \int_0^T \left[\left<\frac{D^2 X}{dt^2}, \frac{D^2 Y}{dt^2} \right> + \left<Y, F(X, \dot{q}) + \nabla_X \grad V \right>\right]dt.
\end{align*}

\begin{definition}
The \text{extended index} of $q$, denoted by $\ind_0(q)$, is the dimension of the largest subspace of $\mathring{H}^2_q$ on which $I$ is negative semidefinite. 
\end{definition}

\begin{lemma}\label{lemma: index}
$\ind_0(q)$ is finite.
\end{lemma}
\textit{Proof:}
Assume towards contradiction that $\ind_0(q)$ is infinite. Then there exists a sequence $(X_k) \subset \mathring{H}^2_q$ such that $I(X_k, X_k) \le 0$, and $(X_k)$ is orthonormal with respect to the $\mathring{H}^1_q$ product given by
$$\left< X, Y \right>_{\mathring{H}^1_q} = \int_0^T \left[ \left<X, Y\right> + \left< \frac{DX}{dt}, \frac{DY}{dt} \right> \right]dt$$
for all $X, Y \in \mathring{H}_q^2$. Since $I(X_k, X_k) \le 0$, we have that
\begin{align*}
    \int_0^T \left\|\frac{D^2 X_k}{dt^2}  \right\|^2 dt &\le \int_0^T \left|\left<X_k, F(X_k, \dot{q}) + \nabla_{X_k} \grad V\right>\right|dt \\
    &\le \int_0^T \|X_k\| \|F(X_k, \dot{q}) + \nabla_{X_k} \grad V\|dt.
\end{align*}

Since $\grad V, \ q,$ and $R$ are smooth, we may consider the family of linear maps $A_t: T_{q(t)} Q \to T_{q(t)} Q$ given by $A_{t}(X) = F(X, \dot{q}) + \nabla_X \grad V(q(t))$, together with the induced norm on $\text{End}(T_{q(t)}Q)$ given by $|A_{t}| = \sup_{\substack{X\in T_{q(t)} Q \\ \|X\| = 1}} \|A_{t}(X)\|$. Since $[0, T]$ is compact, it follows that $\sup_{t \in [0, T]}|A_{t}| < +\infty$. Hence, utilizing the Cauchy-Schwarz inequality,
\begin{align*}
     \int_0^T \|X_k\| \|F(X_k, \dot{q}) + \nabla_{X_k} \grad V\|dt \le \sup_{t \in [0, T]}|A_{t}| \int_0^T \|X_n\|^2 dt \le \sup_{t \in [0, T]}|A_{t}|
\end{align*}

Hence, $(X_k)$ is bounded with respect to the $\mathring{H}^2_q$ norm. By the Rellich–Kondrachov theorem, there exists a subsequence $(X_{k_i})$ which converges in $\mathring{H}^1_q$. However, this is impossible since $(X_{k_i})$ is orthonormal with respect to the product on $\mathring{H}^1_q$, and hence cannot be Cauchy. $\hfill \square$

\begin{corollary}\label{cor: finite_conj}
Suppose that $q$ is a modified cubic polynomial. Then for any $\tau \in [0, T]$, there are a finite number of points which are biconjugate to $\tau$ along $q$.
\end{corollary}
\textit{Proof:} Suppose that there are infinitely many points biconjugate to $\tau$. Then, at least one of the sets $\tau^- = \{t < \tau : \ t, \tau \ \text{biconjugate}\}$ and $\tau^+ = \{t > \tau : \ t, \tau \ \text{biconjugate}\}$ are infinite. Without loss of generality, suppose that $\tau^+$ is infinite, and consider a sequence $(t_k) \subset \tau^+$. Then, for each $k \in \N$, there exists a bi-Jacobi Field $J_k$ along $q$ such that $J_k(\tau) = J_k(t_k) = \frac{D}{dt}J_k(\tau) = \frac{D}{dt}J_k(t_k) = 0$. Consider the sequence of vector fields $X_k \in \mathring{H}^2_q$ defined by $X_k(t) = \begin{cases} J_k(t), &t\in[\tau, t_k] \\
0, &\text{otherwise}
\end{cases}$. Then it is clear that $(X_k) \subset T_q \Omega$ is a linearly independent sequence and $I(X_k, X_k) = 0$ for all $k \in \N$. However this implies that $\ind_0(q)$ is infinite, which contradicts Lemma \ref{lemma: index}. $\hfill \square$

In the following corollary, we show that the bi-exponential map $\biexp_{(p,v)}^{t}$ is not critical provided that $t$ is sufficiently small. This follows from the fact that the set of points biconjugate to any fixed $\tau \in [0, T]$ is finite, and hence consists of isolated points. We use this idea to show that the restrictions of modified cubic polynomials to sufficiently small intervals are unique.

\begin{proposition}\label{cor: loc_uniq}
Suppose that $q$ is a modified cubic polynomial and $\tau \in [0, T]$. Then there exists an interval $[a, b] \subset [0, T]$ containing $\tau$ such that $q\vert_{[a,b]}$ is the unique solution of \eqref{eqq1} on $\Omega_{\xi, \eta}^{[a, b]}$, where $\xi = (q(a), \dot{q}(a)), \ \eta = (q(b), \dot{q}(b)) \in TQ$.
\end{proposition}

\textit{Proof:} It follows easily from Corollary \ref{cor: finite_conj} that there are a finite number of points biconjugate to $t = \tau$. Hence, there exists an interval $[a, b^\ast]$ containing $\tau$ which contains no points biconjugate to $t = \tau$. Similarly, there exists a (possibly smaller) interval $[a, b]$ containing $\tau$ which contains no points biconjugate to $t = a$. Let $p = q(a)$ and $v = \dot{q}(a)$. By Lemma \ref{lemma: conjugate}, $\biexp_{(p, v)}^{|b-a|}$ is not critical at $(\frac{D}{dt}\dot{q}(a), \frac{D^2}{dt^2}\dot{q}(a)) \in (T_p Q)^2$. The result then follows immediately by Lemma \ref{lemma: unique_exp}. $\hfill\square$

Before proving that the modified cubic polynomials are exactly the $Q$-local minimizers of $J$, we need one additional result from \cite{existence}.

\begin{theorem}\label{existence}\cite{existence}
$J$ attains its minimum in $\Omega$.
\end{theorem}

\begin{proposition}\label{prop: loc_min}
$q \in \Omega$ is a modified cubic polynomial if and only if it is a $Q$-local minimizer of $J$. Moreover, $q$ is the locally unique minimizer. 
\end{proposition}

\noindent \textit{Proof:} It is clear that if $q$ is a $Q$-local minimizer of $J$, then it is a modified cubic polynomial. For the other direction, fix $\tau \in [0, T].$ From Proposition \ref{cor: loc_uniq}, there exists a neighborhood $[a, b]$ containing $\tau$ such that $q\vert_{[a,b]}$ is the unique critical point of $J$ on $\Omega_{\xi, \eta}^{[a, b]}$, where $\xi = (q(a), \dot{q}(a)), \ \eta = (q(b), \dot{q}(b)) \in TQ$. However, from Theorem \ref{existence}, there exists a global minimizer of $J$ on $\Omega_{\xi, \eta}^{[a, b]}$—which must itself be a critical point of $J$. Hence, $q\vert_{[a,b]}$ is this minimizer. $\hfill \square$

\section{Conclusion}
Throughout this paper, we derived sufficient conditions for optimality in variationally defined obstacle avoidance problems defined on complete and connection Riemannian manifolds, providing natural extensions to the theory of geodesics and Riemannian cubic polynomials. Local minimizers of the action functional were divided into two classes, the so-called $\Omega$-local minimizers and $Q$-local minimizers, which were subsequently studied individually. The former was understood in the context of bi-Jacobi fields and biconjugate points, and certain robustness and local uniqueness results were obtained. The latter category was shown to be equivalent to the critical points of the action, and some local uniqueness results were similarly obtained. 


\section*{Acknowledgements}
Jacob R. Goodman (\textcolor{blue}{jacob.goodman@icmat.es}) conducts his research at Instituto de Ciencias Matematicas
(CSIC-UAM-UC3M-UCM), Calle Nicolas Cabrera 13-15, 28049, Madrid, Spain. The project that gave rise to these results received the support of a fellowship from ”la Caixa” Foundation (ID 100010434). The fellowship code is LCF/BQ/DI19/11730028. Additionally, support has been given by the ``Severo Ochoa Programme for Centres of Excellence'' in R$\&$D (SEV-2015-0554). The author wishes to express gratitude to Dr. Leonardo Colombo for his invaluable contributions and discussions along the writing of this article. All of the results are original and have not been presented nor submitted to conference.


\end{document}